\documentclass[11pt]{article}

\usepackage[T1]{fontenc}
\usepackage[latin1]{inputenc}

\usepackage[]{graphicx,float,latexsym,times}
\usepackage{amsfonts,amstext,amsmath,amssymb,amsthm,dsfont}
 
\usepackage{enumerate}
\usepackage{setspace}

\RequirePackage{xcolor}
\RequirePackage[numbers]{natbib}
\RequirePackage[colorlinks,citecolor=blue,urlcolor=blue]{hyperref}

\usepackage{setspace}

\usepackage{epsfig}
\usepackage{graphicx}
\usepackage{wrapfig}
\usepackage{amssymb}
\usepackage{amsfonts}
\usepackage{amsmath}
\usepackage{amsthm}
\usepackage{enumerate}
\usepackage{multicol}
\usepackage{bbm} 

\newtheorem{theorem}{Theorem}[section]
\newtheorem{proposition}[theorem]{Proposition}
\newtheorem{lemma}[theorem]{Lemma}

\newtheorem{remark}{Remark}[section]
\newtheorem{corollary}[theorem]{Corollary}
\newcommand\fdem{\hfill $\Box$}
\newcommand\cA{{\cal A}}

\newcommand\cB{{\cal B}}

\newcommand\cP{{\cal P}}

\newcommand\cQ{{\cal Q}}
\newcommand\cS{{\cal S}}
\newcommand\cR{{\cal R}}
\newcommand\cT{{\cal T}}

\newcommand\ve{\varepsilon}

\def\bbr{{\mathbb R}}

\def\bbz{{\mathbb Z}}

\def\text#1{\hbox{#1}}
\def\proof{{\noindent \bf Proof. }}

\def\E{{\bf E}}

\def\c{{\bf c}}
\def\D{{\bf D}}

\def\M{{\bf M}}
\def\m{{\bf m}}
\def\U{{\bf U}}

\def\r{{\bf r}}
\def\l{{\bf l}}

\def\v{{\bf v}}
\def\a{{\bf a}}
\def\d{{\bf d}}
\def\b{{\bf b}}
\def\q{{\bf q}}
\def\p{{\bf p}}
\def\z{{\bf z}}
\def\L{{\bf L}}

\def\x{{\bf x}}
\def\e{{\bf e}}
\newcommand{\wh}{\widehat}
\newcommand{\wt}{\widetilde}

\newcommand\Er{\mbox{Err}}

\newcommand\Rep{\mbox{Re}}
\def\R{{\bf R}}
\def\Chi{{\bf 1}}
\def\d{\mathrm{d}}
\def\build #1_#2{\mathrel{\mathop{\kern 0pt #1}\limits_\zs{#2}}}
\newcommand{\zs}[1]{{\mathchoice{#1}{#1}{\lower.25ex\hbox{$\scriptstyle#1$}}
{\lower0.25ex\hbox{$\scriptscriptstyle#1$}}}}
\numberwithin{equation}{section}

\begin{document}
\bibliographystyle{alpha}

\title{
{\bf Sharp non asymptotic oracle inequalities for non parametric 
computerized tomography  model
}
\thanks{
This work was done under 
 the RSF grant 17-11-01049 (National Research Tomsk State University).
}
}

\author{Fourdrinier, D.
\thanks{
 LITIS, EA 4108, Universit\'e de Rouen, France and
International Laboratory SSP \& QF,
 National Research Tomsk State University, Russia}
and
\smallskip
Pergamenshchikov S.M.\thanks{
 Laboratory of Mathematics LMRS,
  University of Rouen, France and
International Laboratory SSP \& QF,
 National Research Tomsk State University, Russia
 e-mail: Serge.Pergamenchtchikov@univ-rouen.fr}
}


\maketitle

\begin{abstract}
We consider non parametric estimation problem
for stochastic tomography regression model, i.e. we consider the estimation problem of  function of multivariate variables (image)
 observed through its Radon transformation calculated with the random errors. For this problem we develop a new adaptive
 model selection method. By making use the Galtchouk and Pergamenshchikov approach we construct the model selection procedure for which we show a
 sharp non asymptotic oracle inequality for the both 
 usual and robust quadratic risks, i.e. we show that the proposed procedure is optimal in the oracle inequalities sense.
\end{abstract}

\vfill



\eject


\section{Introduction} \label{sec1}

In this paper we consider the multivariate regression model proposed in 
\cite{Natterer1986} for
computerized tomography problems, i.e. we consider the following regression model
\begin{equation}\label{equa1.6}
 y(\nu,\varsigma) = R(S)(\nu,\varsigma) + \xi(\nu,\varsigma) 
\end{equation}
where
$R(\cdot)$ is the Radom transformation, $S$ is a $\bbr^{d}\to\bbr$ function from
$\L_\zs{2}(\bbr^{d})$ such that $S(x)=0$ for $\vert x\vert\ge \x_\zs{*}$ for some fixed $\x_\zs{*}>0$. The vector $\nu\in\bbr^{d}$
with $\vert\nu\vert=1$, $\varsigma\in\bbr$
 and $\xi(\nu,\varsigma)$ is the noise variable. 
Our aim, in this paper, is to estimate the function $S$ based on the vector of 
observations 
\begin{equation}\label{equa1.6++}
y_\zs{l}
=R(S)(\nu_\zs{l},\varsigma_\zs{l})
+\xi_\zs{l}\,,\quad 1\le l\le n\,,
\end{equation}
where $(\xi_\zs{l})_\zs{1\le l\le n}$
is i.i.d. sequence 
with the unknown distribution function $\p$ under which
\begin{equation}\label{noise_conds}
\E_\zs{\p}\,\xi_\zs{1}=0\,,\quad \E_\zs{\p}\,\xi^{2}_\zs{1} = \sigma_\zs{\p}
\quad\mbox{and}\quad
\E_\zs{\p}\,\xi^{4}_\zs{1} <\infty
\,.
\end{equation}

We assume that the distribution $Q$ belongs to some distribution family $\cQ_\zs{n}$ which will be specified below.

 In this case
we use the robust estimation approach 
proposed  in \cite{GaltchoukPergamenshchikov2006, KonevPergamenshchikov2012, KonevPergamenshchikov2015}
for the nonparametric estimation. According to this approach we have to construct an estimator $\wh{S}_\zs{n}$
 (any function of $(y_\zs{l})_\zs{0\le l\le n}$) for $S$ to minimize the robust risk defined as
\begin{equation}\label{sec:risks}
\cR^{*}_\zs{n}(\wh{S}_\zs{n},S)=\sup_\zs{\p\in\cP_\zs{n}}\,
\cR_\zs{\p}(\wh{S}_\zs{n},S)\,,
\end{equation}
where  $\cR_\zs{\p}(\cdot,\cdot)$ is the usual quadratic risk of the form
\begin{equation}\label{sec:risks_00}
\cR_\zs{\p}(\wh{S}_\zs{n},S):=
\E_\zs{\p}\,\|\wh{S}_\zs{n}-S\|^2
\quad\mbox{and}\quad
\|S\|^2
=
 \int_\zs{[-\x_\zs{*},\x_\zs{*}]^d} |S(x)|^2 \d x 
\,.
\end{equation}  
It is clear that if we don't know the distribution
of the observation one  needs to find an estimator which will be optimal for all possible observation distributions. 
Moreover in this paper we consider the estimation problem in the adaptive setting,
i.e. when the regularity of $S$ is unknown. 
To this end
 we use the adaptive method based on the model selection approach. 
The interest to such statistical procedures is explained by the fact that
they provide
adaptive solutions for the nonparametric estimation through
oracle inequalities which give the  non-asymptotic upper
bound for the quadratic risk including
the minimal risk over chosen  family of estimators.
It should be noted that for the first time the model selection methods
were proposed by
 Akaike \cite{Akaike1974} and Mallows \cite{Mallows1973}
for parametric models. Then,
these methods had been developed
for the nonparametric estimation
 and
the oracle inequalities for the quadratic risks was obtained  by
Barron, Birg\'e and Massart \cite{BarronBirgeMassart1999},
by Fourdrinier and Pergamenshchikov \cite{FourdrinierPergamenshchikov2007}
for the regression models in discrete time
and
\cite{KonevPergamenshchikov2010} in continuous time.
Unfortunately, the oracle inequalities obtained in these papers
 can not
provide the efficient estimation in the adaptive setting, since
the upper bounds in these inequalities
have some fixed coefficients in the main terms which are more than one.
To obtain  the efficiency property for estimation
procedures  one has to obtain
 the sharp oracle inequalities, i.e.
 in which
  the factor at the principal term on the right-hand side of the inequality
 is close to unity.
The first result on sharp inequalities is most likely due to Kneip \cite{Kneip1994}
who studied a Gaussian regression model in the discrete time.
It will be observed that the
derivation of oracle inequalities usually rests upon
the fact that the initial model, by applying the Fourier transformation,
can be reduced to the Gaussian independent observations.
However, such  transformation is possible only for
Gaussian models with independent homogeneous observations or for
 inhomogeneous ones with  known correlation characteristics.
For the general non Gaussian observations one needs to use the approach proposed by
 Galtchouk and Pergamenshchikov
 \cite{GaltchoukPergamenshchikov2009a, GaltchoukPergamenshchikov2009b}  for the heteroscedastic
 regression models in discrete time and developed then by
Konev and Pergamenshchikov in
\cite{KonevPergamenshchikov2009a, KonevPergamenshchikov2009b, KonevPergamenshchikov2012, KonevPergamenshchikov2015}
for semimartingale models in  continuous time.
In general the model selection is an
adaptive rule $\wh{\lambda}$ which choses an estimator $S^{*}=\wh{S}_\zs{\wh{\lambda}}$ from 
an estimate family 
$(\wh{S}_\zs{\lambda})_\zs{\lambda\in\Lambda}$. The goal of  this selection is to prove the following
nonasymptotic 
oracle inequality: for any sufficient small $\delta >0$ and any observation duration $n\ge 1$
\begin{equation}\label{sec:In.3+5}
\cR_\zs{\p}( S^{*},S)\,
\le\,(1+\delta)\,\min_\zs{\lambda\in\Lambda}\,\cR_\zs{\p}(\wh{S}_\zs{\lambda},S)+\delta^{-1}
\cB_\zs{n}\,,
\end{equation}
where the rest term
$\cB_\zs{n}$  
 is sufficiently small with respect to the minimax convergence rate.
Such oracle inequalities are called {\sl sharp}, since the coefficient in the main term 
$1+\delta$
is close to one for sufficiently small $\delta>0$. 
The rest of the paper is organized as follows. 
In Section \ref{secMnSlsprs} we state the main conditions for the model \eqref{equa1.6} and we construct the model selection procedures.
In Section \ref{sec:Mrs} we give the main results on the oracle inequalities.
In Section \ref{sec:Prsm} we study the main properties of the model \eqref{equa1.6}. In Section \ref{sec:Pr}
we prove all results. Appendix \ref{sec:A} contains all technical and auxiliary  proofs.

\bigskip

\section{Model selection} \label{secMnSlsprs}
We assume that the noise distribution $\p$ belong to the probability family $\cP_\zs{n}$ is defined as
\begin{equation}
\label{family_prb}
\varsigma_\zs{*}\le \sigma_\zs{\p}\le\varsigma^{*}
\quad\mbox{and}\quad
\E_\zs{\p}\xi^{4}_\zs{1}
\le \varsigma^{*}_\zs{1}
\,,
\end{equation}
where the unknown 
bounds $0<\varsigma_\zs{*}\le \varsigma^{*}$ and $\varsigma^{*}_\zs{1}$ are functions of $x_\zs{*}$, i.e.  
$\varsigma_\zs{*}=\varsigma_\zs{*}(n)$, $\varsigma^{*}=\varsigma^{*}(n)$ and  $\varsigma^{*}_\zs{1}=\varsigma^{*}_\zs{1}(n)$, such that
for any $\check{\epsilon}>0,$
\begin{equation}\label{sec:Mrs.5-1}
\lim_\zs{n\to\infty}n^{\check{\epsilon}}\,\varsigma_\zs{*}(n)=+\infty
\quad\mbox{and}\quad
\lim_\zs{n\to\infty}\,\frac{\varsigma^{*}(n)+\varsigma^{*}_\zs{1}(n)}{n^{\check{\epsilon}}}=0\,.
\end{equation}

\noindent 
First,   we define
the trigonometric basis 
in $\L_\zs{2}[-\x_\zs{*},\x_\zs{*}]$ as 
\begin{equation}\label{equa1.7} 
 \varphi_\zs{k}(\varsigma) 
  = 
 \frac{1}{\sqrt{2 \x_\zs{*}}} \, 
 e^{i k \varsigma\pi/\x_\zs{*}}
 \quad\mbox{and}\quad
i=\sqrt{-1}\,.
\end{equation}
For any $x = (x_\zs{1}, \ldots ,x_\zs{d})\in [-\x_\zs{*},\x_\zs{*}]^{d}$ 
 we can represent  the function $S$ in $\L_\zs{2}[-\x_\zs{*},\x_\zs{*}]$ 
as
\begin{equation}\label{equa1.14}
 S(x) = \sum_{j=(j_1, \ldots ,j_d) \in \bbz^d} \theta_\zs{j} \, \Phi_\zs{j}(x) 
\,,
\end{equation}
where  
$$
 \Phi_\zs{j}(x) 
 = \prod_{l=1}^d \varphi_{j_\zs{l}}(x_\zs{l}) 
\quad\mbox{and}\quad
 \theta_\zs{j} = 
 \int_{[-\x_\zs{*},\x_\zs{*}]^d} S(z) \, \Phi_\zs{j}(z) \, \d z \,. 
$$
So, taking into account the equality
\eqref{equa1.4} we can represent the coefficient $\theta_\zs{j}$
 as
\begin{equation*}\label{equa1.14++}
\theta_\zs{j}=\frac{1}{(\sqrt{2\x_\zs{*}})^{d}}\,\int^{\x_\zs{*}}_\zs{-\x_\zs{*}}\,R(S)(\nu_\zs{j},\varsigma)e^{i\beta_\zs{j} \varsigma}\d \varsigma
\,,
\end{equation*}
where $\beta_\zs{j}=\vert j\vert \pi/\x_\zs{*}$ and $\nu_\zs{j}=j/\vert j\vert$ for $\vert j\vert>0$.
Taking into account here the property \eqref{FinRadonTr_1}
 we obtain that for any $L_\zs{j}\ge x_\zs{*}$
the Fourier coefficient can be rewritten as
\begin{equation}\label{equa1.14++15}
\theta_\zs{j}
=
\frac{1}{(\sqrt{2\x_\zs{*}})^{d}}\,
\int^{L_\zs{j}}_\zs{-L_\zs{j}}\,R(S)(\nu_\zs{j},\varsigma)e^{i\beta_\zs{j} \varsigma}\d \varsigma
\,,
\end{equation}
\noindent
To calculate this coefficient we use the approximation 
\begin{equation}\label{equa2.6}
 \a_\zs{j}
  =
\frac{1}{(\sqrt{2\x_\zs{*}})^{d}}\,
 \sum_{l=1}^{\q} R(S)(\nu_\zs{j},s_\zs{l}) \,\psi_\zs{j,l}
\quad\mbox{and}\quad
 \psi_\zs{j,l}=
  \int_{s_\zs{j,l-1}}^{s_\zs{j,l}} 
  e^{i\beta_\zs{j}\varsigma}
   \d \varsigma
\,,
\end{equation}
where 
 $(s_\zs{l})_\zs{1\le l\le n}$ is the uniform partition of the interval $[-L_\zs{j}\,,\,L_\zs{j}]$, i.e. 
$$
s_\zs{j,l} = -L_\zs{j} + \frac{2 L_\zs{j}}{\q}l\,.
$$
The number of points $\q=\q_\zs{n}$ will be chosen later.
Using this approximation we estimate $\theta_\zs{j}$ as
\begin{equation}\label{equa2.7}
 \wh{\theta}_\zs{j} = \sum_\zs{l=1}^{\q} y_\zs{j,l}
 \psi_\zs{k,l}
 \quad\mbox{and}\quad
  y_\zs{j,l}=y(\nu_\zs{j},s_\zs{j,l})
  \,.
\end{equation}
Therefore, it can be represented as 
\begin{equation}\label{equa2.8}
 \wh{\theta}_\zs{j}
 =\theta_\zs{j} + \zeta_\zs{j}
 \quad\mbox{and}\quad
 \zeta_\zs{j}=\b_\zs{j}+\frac{1}{\sqrt{\q}}\,\eta_\zs{j}
 \,,
\end{equation}
where $\b_\zs{j}=\a_\zs{j}-\theta_\zs{j}$,
\begin{equation}\label{equa2.11}
\eta_\zs{j} = \sqrt{\q} \, 
 \sum_\zs{l=1}^{\q} \xi_\zs{j,l}
 \psi_\zs{j,l}
 \quad\mbox{and}\quad
\xi_\zs{j,l}=
 \xi(\nu_\zs{j},s_\zs{j,l})
 \,. 
\end{equation}
Note that as it shown in Proposition
\ref{Pr.sec:Prsm.1} the second moment of this random variable is given as
\begin{equation}\label{variance_eta}
\E_\zs{\p}\,\vert \eta_\zs{j}\vert^{2}=
\sigma_\zs{\p}
\varpi_\zs{j}
\,,
\end{equation}
where 
$$
\varpi_\zs{j}
=4L^{2}_\zs{j}\,
\frac{\sin^{2}\check{\beta}_\zs{j}}{\check{\beta}^{2}_\zs{j}}
\quad\mbox{and}\quad
\check{\beta}_\zs{j}=\pi
(1+[1/\nu^{*}_\zs{j}])\nu^{*}_\zs{j}\,
\frac{\vert j\vert}{\q}\,.
$$
To obtain the uncorrelated property for the random variables $\eta_\zs{j}$
  we set $L_\zs{j}$ as
\begin{equation}  \label{L-jdfen}
L_\zs{j}=L(\nu_\zs{j})=(1+[1/\nu^{*}_\zs{j}])\nu^{*}_\zs{j} \x_\zs{*}\,,
\end{equation}  
  where the coefficient $\nu^{*}_\zs{j}$ is the absolute value of the first nonzero component of the vector 
  $\nu_\zs{j}=(\nu_\zs{j,1},\ldots,\nu_\zs{j,d})'$, i.e. $\nu^{*}_\zs{j}=\min\{\vert \nu_\zs{j,l}\vert >0\,,1\le l\le d\}$.
Taking into account that in this case $L_\zs{j}\ge \x_\zs{*}$ we obtain that
\begin{equation}
\label{infvarpi--}
\varpi_\zs{*}
=
\inf_\zs{\vert j\vert\le \q}\,\varpi^{2}_\zs{j}>0\,.
\end{equation}

So, using the estimators \eqref{equa2.7}
we will estimate the function $S$. The idea is the following, first we replace the infinite sum in \eqref{equa1.14}
by the finite sum over the set 
\begin{equation}\label{set_summing}
\cS_\zs{n}=\{-\m_\zs{n},\ldots,\m_\zs{n}\}^d\,,
\end{equation}
where the integer $\m_\zs{n}\ge 1$ will be specify below. Then, 
according to the Pinsker weighted least square method 
we will replace the Fourier coefficients in \eqref{equa1.14}
with its estimators \eqref{equa2.7} multiplied by some coefficient $0\le \lambda(j)\le 1$, i.e.

\begin{equation}\label{sec:Or.7}
\wh{S}_\zs{\lambda}(x)=\sum_{j \in \cS_\zs{n}} \lambda(j) \, \widehat{\theta}_\zs{j} \, 
\Phi_\zs{j}(x)\,,
\quad x \in [-\x_\zs{*},\x_\zs{*}]^d\,,
\end{equation}
the weight vector $\lambda=(\lambda(j))_\zs{j \in \cS_\zs{n}}$
belongs to some finite set $\Lambda$ from $[0,1]^{r_n}$ and $r_\zs{n}=(2\m_\zs{n}+1)^{d}$. 
We set
 \begin{equation}\label{sec:Mo.2}
\check{\iota}=\mbox{card}(\Lambda)
\quad\mbox{and}\quad
\vert\Lambda\vert_\zs{*}=\max_\zs{\lambda\in\Lambda}\,
\left(\check{L}(\lambda)
+
\check{L}_\zs{1}(\lambda)
\right)
\,,
\end{equation}
where $\check{L}(\lambda)=\sum_\zs{j\in\bbz^{d}}\lambda(j)\varpi_\zs{j}$
and 
$\check{L}_\zs{1}(\lambda)=\sum_\zs{j\in\bbz^{d}}\Chi_\zs{\{\lambda\neq 0\}}$.

Now we need to write a cost function to choose a weight $\lambda\in\Lambda$. Of course,
it is obvious, that the best way is to minimize the cost function which is equal to 
the empirical squared error
$$
\Er_\zs{n}(\lambda)=\|\widehat{S}_\zs{\lambda}-S\|^2 
 = \int_{[-\x_\zs{*},\x_\zs{*}]^d} |\widehat{S}_\zs{\lambda}(x) - S(x)|^2 \d x \,,
$$
 which in our case is equal to
\begin{equation}\label{sec:Or.8}
\Er_\zs{n}(\lambda) = 
\sum_{j \in \bbz^{d}} \lambda^2(j) \, \vert\widehat{\theta}_\zs{j}\vert^{2}\,
-
2\,\sum_{j \in \bbz^{d}} \lambda(j)
\,\Rep\,
\bar{\theta}_\zs{j}\,
\widehat{\theta}_\zs{j}\, + 
\Vert S\Vert^{2}\,,
\end{equation}
where $\lambda(j)=0$ for $j\in\bbz^{d}\setminus \cS_\zs{n}$.
Since the coefficients $\theta_\zs{j}$ are unknown, we need to replace the term
$\widehat{\theta}_\zs{j}\,\bar{\theta}_\zs{j}$ by an estimate which we choose as
\begin{equation}
\label{coeff_sigma}
\wt{\theta}_\zs{j}=
\vert \widehat{\theta}_\zs{j}\vert^{2}-\frac{\wh{\sigma}}{\q} \, \varpi_\zs{j}\,.
\end{equation}
Here $\wh{\sigma}$ is an estimate for $\sigma_\zs{\p}$ which is given in 
\eqref{sec:wh_sigma}. If the variance $\sigma_\zs{\p}$ is known we set 
$\wh{\sigma}=\sigma_\zs{\p}$, otherwise, we can choose as
\begin{equation}\label{sec:wh_sigma}
\wh{\sigma}=
\frac{1}{\wt{\q}}
\,
\sum_\zs{j\in \cT_\zs{n}}\,
\vert\wh{\theta}_\zs{j}\vert^{2}
\quad\mbox{and}\quad
\wt{\q}=\frac{\sum_\zs{j\in \cT_\zs{n}}\varpi_\zs{j}}{\q}
\,,
\end{equation}
where $\cT_\zs{n}=\{[\sqrt{\m_\zs{n}}]+1,\ldots,\m_\zs{n}\}^{d}$
and the coefficients $\varpi_\zs{j}$ are defined in
\eqref{coeff_sigma}.

Moreover, for this substitution to the empirical squared error one needs to pay
a penalty. Finally, we define the cost function in the following way
\begin{equation}\label{sec:Or.9}
J_\zs{n}(\lambda) = \sum_{j \in \bbz^{d}} \lambda^2(j) \, \vert\wh{\theta}_\zs{j} \vert^{2}-
2 \sum_{j \in \bbz^{d}}\,\lambda(j)\,\wt{\theta}_\zs{j}
+ \delta \, \wh{P}_\zs{n}(\lambda)\,,
\end{equation}
where $\rho$ is some positive penalty coefficient which will be chosen later and  the penalty term 
$\wh{P}_\zs{n}(\lambda)$ we choose as
\begin{equation}\label{sec:Or.10}
\wh{P}_\zs{n}(\lambda)=\frac{\wh{\sigma}}{\q}\,\check{L}(\lambda^{2})
\,,
\end{equation}
where $\lambda^{2}=(\lambda^{2}(j))_\zs{j\in\bbz^{d}}$.
In the case when the $\sigma_\zs{\p}$ is known we set
\begin{equation}\label{sec:Or.10++1}
P_\zs{n}(\lambda)=\frac{\sigma_\zs{\p}}{\q}\,\check{L}(\lambda^{2})
\,.
\end{equation}

\noindent 
Now, we define the model selection procedure as
\begin{equation}\label{sec:Mo.9}
\wh{S}_\zs{*} = \wh{S}_\zs{\hat \lambda}
\quad\mbox{and}\quad
\wh{\lambda}= \mbox{argmin}_\zs{\lambda\in\Lambda} J_\zs{n}(\lambda)\,.
\end{equation}
We recall that the set $\Lambda$ is finite so $\hat \lambda$ exists. In the case when $\hat \lambda$ is not unique, we take one of them.

Let us now specify the weight coefficients $\lambda=(\lambda(j))_\zs{j \in \cS_\zs{n}}$. Consider, for some fixed $0<\varepsilon<1,$
a numerical grid of the form
\begin{equation}\label{sec:Ga.0}
\cA=\{1,\ldots,k^*\}\times\{\varepsilon,2\varepsilon,\ldots,[1/\ve^2]\varepsilon\}\,,
\end{equation}
where $[a]$ means the integer part of the number $a$. We assume that both parameters $k^*\ge 1$ and $\varepsilon$ are functions of $x_\zs{*}$, i.e.
$k^*=k^*(n)$ and $\ve=\ve(n)$, such that
\begin{equation}\label{sec:Ga.1}
\left\{
\begin{array}{ll}
&\lim_\zs{n\to\infty}\,k^*(n)=+\infty\,,
\quad\lim_\zs{n\to\infty}\,\dfrac{k^*(n)}{\ln n}=0\,,\\[6mm]
&
\lim_\zs{n\to\infty}\,\varepsilon(n)=0
\quad\mbox{and}\quad
\lim_\zs{n\to\infty}\,n^{\check{\delta}}\ve(n)\,=+\infty
\end{array}
\right.
\end{equation}
for any $\check{\delta}>0$. One can take, for example, for $n\ge 2$
\begin{equation}\label{sec:Ga.1-00}
\ve(n)=\frac{1}{ \ln n }
\quad\mbox{and}\quad
k^*(n)=k^{*}_\zs{0}+\sqrt{\ln n}\,,
\end{equation}
where $k^{*}_\zs{0}\ge 0$ is some fixed constant.
 For each $\alpha=(\beta, \l)\in\cA$, we introduce the weight
sequence
\begin{equation}\label{sec:Ga.2}
\lambda_\zs{\alpha}=(\lambda_\zs{\alpha}(j))_\zs{j \in \cS_\zs{n}}
\quad\mbox{and}\quad
\lambda_\zs{\alpha}(j)=\prod^{d}_\zs{l=1}\,\check{\lambda}_\zs{\alpha}(j_\zs{l})
\end{equation}
with the elements 
$$
\check{\lambda}_\zs{\alpha}(t)=\Chi_\zs{\{1\le t<t_\zs{*}\}}+
\left(1-(t/\omega_\alpha)^\beta\right)\,
\Chi_\zs{\{ t_\zs{*}\le t\le \omega_\zs{\alpha}\}},
$$
where
$t_\zs{*}=1+\left[\ln\upsilon_\zs{n}\right]$, $\omega_\zs{\alpha}=(\check{\d}_\zs{\beta}\,\l\upsilon_\zs{n})^{1/(2\beta+d)}$,
$$
\check{\d}_\zs{\beta}=\frac{(\beta+1)(2\beta+1)}{\pi^{2\beta}\beta}
\quad\mbox{and}\quad
\upsilon_\zs{n}=\m_\zs{n}/\varsigma^{*}
\,.
$$
and the threshold $\varsigma^{*}$ is introduced in \eqref{family_prb}.
Now we define the set $\Lambda$ 
as
\begin{equation}\label{sec:Ga.3}
\Lambda\,=\,\{\lambda_\zs{\alpha}\,,\,\alpha\in\cA\}\,.
\end{equation}
It will be noted that in this case the cardinal of the set $\Lambda$ is  
\begin{equation}
\label{sec:Ga.1++1--1}
\check{\iota}_\zs{n}=k^{*} [1/\ve^2]\,.
\end{equation}
Moreover,
taking into account that $\check{\d}_\zs{\beta}<1$ for $\beta\ge 1$ 
we obtain for the set \eqref{sec:Ga.3}
\begin{equation}
\label{sec:Ga.1++1--2}
 \vert \Lambda\vert_\zs{*}\,
 \le\,
\sup_\zs{\alpha\in\cA}
  (\omega_\zs{\alpha})^{d}
\le (\upsilon_\zs{n}/\ve )^{d/(2+d)}\,.
\end{equation}

\begin{remark} \label{Re.sec.Ex.2_2}
Note that the form \eqref{sec:Ga.2}
for the weight coefficients 
 was proposed by Pinsker in \cite{Pinsker1981}
 for the efficient estimation in the nonadaptive case, i.e. when the regularity parameters of the function $S$ are known. 
In the adaptive case  these weight coefficients are
  used in \cite{KonevPergamenshchikov2012, KonevPergamenshchikov2015}
   to show the asymptotic efficiency for model selection procedures.
\end{remark}

\section{Main results}\label{sec:Mrs}

Now we formulate all non asymptotic oracle inequalities. Before, let us first introduce
the following auxiliary function which is used fto describe the rest terms in the oracle inequalities.
\begin{equation}
\label{psi-qq-p}
\Psi_\zs{\p}=\left(
1+\sigma_\zs{\p}
+
\frac{1}{ \sigma_\zs{\p}} 
\right)
\,\E_\zs{\p}\xi^{4}_\zs{1}
\,
\check{\iota}
\,.
\end{equation}
 First, we obtain the oracle inequality for the risk \eqref{sec:risks_00}.

\begin{theorem}\label{Th.sec:OI.1} 
There exists some  constant $\check{\upsilon}>0$ such that for any
 $0 <\delta< 1/8$, any $\q\ge 2d \m_\zs{n}+2$ and any $\p\in\cP_\zs{n}$  the estimator of $S$ given in \eqref{sec:Mo.9} satisfies the following oracle inequality
\begin{align}\nonumber
 \cR_\zs{\p}(\wh{S}_*,S)&  \le \frac{(1+2\delta) }{1-4\delta} \min_\zs{\lambda\in\Lambda}\cR_\zs{\p}(\wh{S}_\zs{\lambda},S)\\[2mm]\label{or-in-p-00}
&+  \check{\upsilon}
\frac{\Psi_\zs{\p}+\vert\Lambda\vert_\zs{*}\,(\E_\zs{\p}\,|\wh{\sigma} -\sigma_\zs{\p}\vert+\wt{\b}/\q)}{\q\delta}
 \,,
\end{align}
where $\wt{\b}=\q^{2}\sup_\zs{j\in\bbz^{d}}\vert\b_\zs{j}\vert^{2}$.
\end{theorem}
Note that, if $\sigma_\zs{\p}$ is known we obtain the following results.
\begin{corollary}\label{Co.sec:OI.1} 
There exists some  constant $\check{\upsilon}>0$ such that for any
 $0 <\delta< 1/8$, any $\q\ge 2d \m_\zs{n}+2$ and any $\p\in\cP_\zs{n}$  the estimator of $S$ given in \eqref{sec:Mo.9} satisfies the following oracle inequality
\begin{equation}\label{cr-or-in-p-00}
 \cR_\zs{\p}(\wh{S}_*,S)  \le \frac{(1+2\delta) }{1-4\delta}\min_\zs{\lambda\in\Lambda} \cR_\zs{\p}(\wh{S}_\zs{\lambda},S)
+  \check{\upsilon}
\frac{\Psi_\zs{\p}+\vert\Lambda\vert_\zs{*}\wt{\b}/\q}{\q\delta}
 \,.
\end{equation}
\end{corollary}

\noindent Now we  study the estimate \eqref{sec:wh_sigma}.
\begin{proposition}\label{Pr.sec:Si.1}
Assume that  the partial derivative 
$\partial^{d}/\partial x_\zs{1}\dots\partial x_\zs{d}\, S$
is continuous and $\check{\r}$ be the Lipshits constant for $S$. Then, there exists some positive $\check{\upsilon}>0$ such that
for any $\m_\zs{n}\ge 4$, $\q_\zs{*}\ge 1$ and $\m^{d}_\zs{n}\le \q\le \q_\zs{*} \m^{d}_\zs{n}$,
\begin{equation}\label{sec:Si.3}
\E_\zs{\p}|\wh{\sigma}-\sigma_\zs{\p}|
\le
\check{\upsilon}\q_\zs{*}\,
\frac{
(1+\check{\r}^{2}+\check{\tau}_\zs{d}(S))(1+\E_\zs{\p}\xi^{4}_\zs{1})}{\m^{d/2}_\zs{n}}
\,,
\end{equation}
where and $\check{\tau}_\zs{d}(S)=\Vert \partial^{d}/\partial x_\zs{1}\dots\partial x_\zs{d}\, S\Vert^{2}$.
\end{proposition}

\bigskip

\noindent 
Theorem \ref{Th.sec:OI.1} and Proposition \ref{Pr.sec:Si.1}
implies the following result.

\begin{theorem}\label{Th.sec:OI.1_22} 
Assume that  the partial derivative 
$\partial^{d}/\partial x_\zs{1}\dots\partial x_\zs{d}\, S$
is continuous and $\check{\r}$ be the Lipshits constant for $S$. Then, there exists some positive $\check{\upsilon}>0$ such that
for any $\m_\zs{n}\ge 4$, $\q_\zs{*}\ge 1$, $\m^{d}_\zs{n}\le \q\le \q_\zs{*} \m^{d}_\zs{n}$, $\vert\Lambda\vert_\zs{*}\le \m^{d/2}_\zs{n}$
and any $\p\in\cP_\zs{n}$  the estimator of $S$ given in \eqref{sec:Mo.9} satisfies the following oracle inequality
\begin{align}\nonumber
 \cR_\zs{\p}(\wh{S}_*,S)&  \le \frac{(1+2\delta) }{1-4\delta} \min_\zs{\lambda\in\Lambda}\cR_\zs{\p}(\wh{S}_\zs{\lambda},S)\\[2mm]\label{or-in-p-00++}
&+  \check{\upsilon}\,\q_\zs{*}\,
\frac{\Psi_\zs{\p}(1+\check{\r}^{2}+\check{\tau}_\zs{d}(S))}{\q\delta}
 \,.
\end{align}
\end{theorem}

\bigskip

The next result presents the non-asymptotic oracle inequality for the robust risk \eqref{sec:risks}
 for the model selection procedure \eqref{sec:Mo.9}, considered with the coefficients
\eqref{sec:Ga.2}.  Using the definition of the probability family $\cP_\zs{n}$ in \eqref{family_prb}
and the function \eqref{psi-qq-p} we can obtain directly the following result.

\begin{theorem}\label{Th.sec:OI.1_22*} 
Assume that  the partial derivative 
$\partial^{d}/\partial x_\zs{1}\dots\partial x_\zs{d}\, S$
is continuous and $\check{\r}$ be the Lipshits constant for $S$. Then, there exists some positive $\check{\upsilon}>0$ such that
for any $\m_\zs{n}\ge 4$, $\q_\zs{*}\ge 1$, $\m^{d}_\zs{n}\le \q\le \q_\zs{*} \m^{d}_\zs{n}$,   the estimator of $S$ given in \eqref{sec:Mo.9} satisfies the following oracle inequality
\begin{align}\nonumber
 \cR^{*}_\zs{n}(\wh{S}_*,S)&  \le \frac{(1+2\delta) }{1-4\delta} \min_\zs{\lambda\in\Lambda}\cR^{*}_\zs{n}(\wh{S}_\zs{\lambda},S)\\[2mm]\label{or-in-p-00++*}
&+  \check{\upsilon}\,\q_\zs{*}\,
\frac{\Psi^{*}_\zs{n}(1+\check{\r}^{2}+\check{\tau}_\zs{d}(S))}{\q\delta}
 \,,
\end{align}
where the coefficient $\Psi^{*}_\zs{n}>0$ is such that for any $\check{\delta}>0$,
\begin{equation}\label{sec:Mrs.7-25.3}
\lim_\zs{n\to\infty}
\,
\frac{\Psi^{*}_\zs{n}}{n^{\check{\delta}}}
=0\,.
\end{equation}
\end{theorem}

\begin{remark}\label{Re.3.1}
Note that the principal term in the right-hand side of the inequality \eqref{or-in-p-00++*} is best
in the class of estimators $(\wh{S}_\zs{\lambda}\,,\,\lambda\in\Lambda)$.
Inequalities of such type are called {\rm the sharp non-asymptotic
oracle inequalities}. The inequality is sharp in the sense that the coefficient
of the principal term may be chosen as close to $1$ as desired.  
\end{remark}

\section{Properties of the regression model \eqref{equa1.6}}
\label{sec:Prsm}


Firs we study the property of the random variables \eqref{equa2.11}

\begin{proposition}
\label{Pr.sec:Prsm.1}
For any vectors $j$ and $k$ from $\bbz^{d}$
\begin{equation}\label{crl_11_0}
\E_\zs{\p}\,\eta_\zs{j}\,\bar{\eta}_\zs{k}=
\sigma_\zs{\p}
\varpi_\zs{j}
\Chi_\zs{\{j=k\}}
\,,
\end{equation}
where $\varpi_\zs{j}$ is defined \eqref{variance_eta}.
 \end{proposition}
 \proof 
  Note that, in view of  
\eqref{equa2.11} we obtain that
$$
\E_\zs{\p}\,\eta_\zs{j}\,\bar{\eta}_\zs{k}
=\q
\sum_\zs{l,l_\zs{1}=1}^{\q}
 \psi_\zs{j,l}
 \psi_\zs{k,l_\zs{1}}
 \E_\zs{\p}\,\xi_\zs{j,l}
\xi_\zs{k,l_\zs{1}}
\,.
$$
It is clear that, if $\nu_\zs{j}\neq \nu_\zs{k}$, then $\E_\zs{\p}\,\xi_\zs{j,l}
\xi_\zs{k,l_\zs{1}}=0$ and, therefore, $\E_\zs{\p}\,\eta_\zs{j}\,\bar{\eta}_\zs{k}
=0$.
Let now $\nu_\zs{j}= \nu_\zs{k}$, but $j\neq k$. In this case
we obtain that
\begin{equation}
\label{crSxi_check}
\E_\zs{\p}\,\eta_\zs{j}\,\bar{\eta}_\zs{k}
=\q\sigma_\zs{\p}\,
\sum_\zs{l=1}^{\q}
 \psi_\zs{j,l}
 \bar{\psi}_\zs{k,l}
\,.
\end{equation}
The functions $\psi_\zs{j,l}$ can be represented as
\begin{equation}
\label{varphi_int}
\psi_\zs{j,l}=e^{i\beta_\zs{j}s_\zs{j,l}}
\,\Upsilon_\zs{j}(\Delta_\zs{j})
\,,
\end{equation}
where $\Upsilon_\zs{j}(z)=\int^{z}_\zs{0}\,e^{-i\beta_\zs{j}x}\,\d x$ and 
$\Delta_\zs{j}=2L_\zs{j}/\q$. Note now, that in view of 
\eqref{L-jdfen} we get that $L_\zs{j}=L_\zs{k}$ for $\nu_\zs{j}= \nu_\zs{k}$. So, in this case for $j\neq k$
\begin{align*}
\sum_\zs{l=1}^{\q}
 \psi_\zs{j,l}
 \bar{\psi}_\zs{k,l}
&=\Upsilon_\zs{j}(\Delta_\zs{j})\Upsilon_\zs{k}(\Delta_\zs{j})
\sum^{\q}_\zs{l=1}\,e^{i(\beta_\zs{j}-\beta_\zs{k})s_\zs{j,l}}
\\[2mm]
&=
\Upsilon_\zs{j}(\Delta_\zs{j})\Upsilon_\zs{k}(\Delta_\zs{j})\,Q_\zs{j}\,
\frac{1-Q^{\q}_\zs{j}}{1-Q_\zs{j}}
\,e^{-i(\beta_\zs{j}-\beta_\zs{k})L_\zs{j}}
\,,
\end{align*}
where $Q_\zs{j}=e^{i(\beta_\zs{j}-\beta_\zs{k})\Delta_\zs{j}}$.  Taking into account the definition of $\nu^{*}_\zs{j}$ in \eqref{L-jdfen},
we obtain that for $\nu_\zs{j}=\nu_\zs{k}$ the difference 
$$
l^{*}=
(\vert j\vert- \vert k\vert)\nu^{*}_\zs{j}
= \vert j\vert\nu^{*}_\zs{j}
- \vert k\vert\nu^{*}_\zs{k}
\in\bbz
\,.
$$
Therefore, taking into account again the definition of $L_\zs{j}$ in in \eqref{L-jdfen}, we obtain that
$$
Q^{\q}
=
e^{
i(\beta_\zs{j}-\beta_\zs{k})\Delta_\zs{j}q}
=e^{i2\pi l^{*}}
=1\,,
$$
i.e. 
\begin{equation}
\label{prd_scale}
\sum_\zs{l=1}^{\q}
 \psi_\zs{j,l}
 \bar{\psi}_\zs{k,l}=\,
 \q\,
 \vert\Upsilon_\zs{j}(\Delta_\zs{j})\vert^{2}
 \Chi_\zs{\{k=j\}}
 \end{equation}
Therefore,
$$
\E\vert \eta_\zs{j}\vert^{2}=
\sigma_\zs{\p}
\varpi_\zs{j}
\quad\mbox{and}\quad
\varpi_\zs{j}=\q^{2}\,\vert\Upsilon_\zs{j}(\Delta_\zs{j})\vert^{2}
\,.
$$
This implies \eqref{crl_11_0}. 
Hence Proposition \ref{Pr.sec:Prsm.1}.
\qed

\begin{proposition}
\label{Pr.sec:Prsm.2}
Let $\z=(\z_\zs{j})_\zs{j\in\bbz^{d}}$ be a family of non random complex numbers. 
Then
\begin{equation}\label{crl_11_0++}
\E_\zs{\p}\,\vert \sum_\zs{j\in\bbz^{d}}\z_\zs{j}\eta_\zs{j}\vert^{2}
\le
\check{\varpi}_\zs{\p}\,\vert \z\vert^{2}
\,,
\end{equation}
where $\check{\varpi}_\zs{\p}=16\sigma_\zs{\p} \x_\zs{*}$ and $\vert \z\vert^{2}=\sum_\zs{j\in\bbz^{d}}\,\vert\z_\zs{j}\vert^{2}$. 
\end{proposition}

\noindent
Now for any non random family of real numbers $x=(x_\zs{j})_\zs{j\in\bbz^{d}}$ with $\vert x\vert<\infty$ we set
\begin{equation}\label{func_U_++}
\U(x)=
\sum_\zs{j\in\bbz^{d}}
\,x_\zs{j}
\,
\wt{\eta}_\zs{j}
\,,
\end{equation}
where $\wt{\eta}_\zs{j}=\vert \eta_\zs{j} \vert^{2}-\sigma_\zs{\p}\varpi_\zs{j}$.

\begin{proposition}
\label{Pr.sec:Prsm.3}
The function \eqref{func_U_++} admits the following upper bound
\begin{equation}\label{crl_12_0++}
\sup_\zs{1\le \#(x)\le \q}
\frac{
\E_\zs{\p}\,\vert 
\U(x)
\vert^{2}
}
{
\vert x\vert^{2}
}
\le
\c_\zs{*}\,
\,\E_\zs{\p}\xi^{4}_\zs{1}
\,,
\end{equation}
where $\#(x)=\sum_\zs{j}\Chi_\zs{\{x_\zs{j}\neq 0\}}$ and  $\c_\zs{*}=5\x_\zs{*}^{2}2^{9}$. 
\end{proposition}
\proof
First, note that  
the random variable $\wt{\eta}_\zs{j}$ can be represented as
$$
\wt{\eta}_\zs{j}=\q\sum^{\q}_\zs{l=1}\,\wt{\xi}_\zs{j,l}
\vert \psi_\zs{j,l}\vert^{2}
+2\q\sum^{\q}_\zs{l=2}\v_\zs{j,l}\xi_\zs{j,l}
\,,
$$
where $\wt{\xi}_\zs{j,l}=\xi^{2}_\zs{j,l}-\sigma_\zs{\p}$ and $\v_\zs{j,l}=\psi_\zs{j,l}\sum^{l-1}_\zs{t=1}\bar{\psi}_\zs{j,t}\xi_\zs{j,t}$.
Moreover, in view of the equallity \eqref{varphi_int}
we obtain that
$$
\wt{\eta}_\zs{j}
=\q\vert\Upsilon_\zs{j}(\Delta_\zs{j})\vert^{2}
\sum^{\q}_\zs{l=1}\,\wt{\xi}_\zs{j,l}
+2\q\sum^{\q}_\zs{l=2}\v_\zs{j,l}\xi_\zs{j,l}
\,.
$$
Therefore, setting
$$
\U_\zs{1,l}(x)=\sum_\zs{j}\,x_\zs{j}\,\vert\Upsilon_\zs{j}(\Delta_\zs{j})\vert^{2}\wt{\xi}_\zs{j,l}
\quad\mbox{and}\quad
\U_\zs{2,l}(x)=\sum_\zs{j}\,x_\zs{j}\,\v_\zs{j,l}\xi_\zs{j,l}
\,,
$$
 we can represent the function \eqref{func_U_++} as
\begin{equation}
\label{repU_eq}
\U(x) =\q \U_\zs{1}(x)+
2\q \U_\zs{2}(x)\,,
\end{equation}
where
$
\U_\zs{1}(x)=\sum^{\q}_\zs{l=1}\,\U_\zs{1,l}(x)$ and
$\U_\zs{2}(x)=\sum^{\q}_\zs{l=1}\,\U_\zs{2,l}(x)$.
Taking into account that the random variables $(\wt{\xi}_\zs{j,l})_\zs{1\le l\le \q}$ are independent with $\E\wt{\xi}_\zs{j,l}=0$, we obtain that
$$
\E_\zs{\p}\,\U^{2}_\zs{1}(x)=
\sum^{\q}_\zs{l=1}\,\E_\zs{\p}\,\U^{2}_\zs{1,l}(x)
\quad\mbox{and}\quad
\E_\zs{\p}\,\vert \U_\zs{2}(x)\vert^{2}=
\sum^{\q}_\zs{l=1}\,\E_\zs{\p}\,\vert\U_\zs{2,l}(x)\vert^{2}
\,.
$$
Now, using here that $\vert\Upsilon_\zs{j}(\Delta_\zs{j})\vert\le \Delta_\zs{j}$, we obtain through the
Cauchy-Bunyakovsky-Schwarz inequality that
$$
\E_\zs{\p}\,\U^{2}_\zs{1,l}(x)\le \Delta^{4}_\zs{j}
\vert x\vert^{2}
\#(x)\,\E_\zs{\p}\,\xi^{4}_\zs{1}
\le 
\,
\frac{2^{8}\x_\zs{*}^{2}\E_\zs{\p}\,\xi^{4}_\zs{1}}{\q^{4}}
\vert x\vert^{2}
\#(x)
\,.
$$
Therefore, for $\#(x)\le \q$ we obtain that
$$
\E_\zs{\p}\,\U^{2}_\zs{1}(x)\le 2^{8}\x_\zs{*}^{2}\E_\zs{\p}\,\xi^{4}_\zs{1}\vert x\vert^{2}\q^{-2}\,.
$$

Moreover, to estimate the last term in \eqref{repU_eq} note that $\U_\zs{2,l}(x)$ can be rewritten as
$$
\U_\zs{2,l}(x)=\sum^{l-1}_\zs{t=1}
\tau_\zs{t,l}
\quad\mbox{and}\quad
\tau_\zs{t,l}=\sum_\zs{j}\,x_\zs{j}\,
\psi_\zs{j,l}\,
\bar{\psi}_\zs{j,t}\xi_\zs{j,t}\,
\xi_\zs{j,l}
\,.
$$ 
It is easy to see that for any $1\le t, s\le l-1$
$$
\E_\zs{\p}\,\tau_\zs{t,l}\bar{\tau}_\zs{s,l}=
\sum_\zs{j,k}\,x_\zs{j}\,x_\zs{k}\,
\psi_\zs{j,l}\,
\bar{\psi}_\zs{j,t}\,
\psi_\zs{k,s}\,
\bar{\psi}_\zs{k,l}\,
\E_\zs{\p}\,\xi_\zs{j,t}\,
\xi_\zs{k,s}\,
\E_\zs{\p}\,
\xi_\zs{j,l}\,
\xi_\zs{k,l}\,,
$$
i.e. $\E_\zs{\p}\,\tau_\zs{t,l}\bar{\tau}_\zs{s,l}=0$ for $t\neq s$ and
$$
\E_\zs{\p}\,\vert\tau_\zs{t,l}\vert^{2}=\sigma^{4}_\zs{\p}
\sum_\zs{j,k}\,x_\zs{j}\,x_\zs{k}\,
\bar{\psi}_\zs{k,l}\,
\psi_\zs{j,l}\,
\bar{\psi}_\zs{j,t}\,
\psi_\zs{k,t}\,
\Chi_\zs{\{\nu_\zs{j}=\nu_\zs{k}\}}
\,.
$$
Therefore, 
\begin{align*}
\E_\zs{\p}\,\vert\U_\zs{2,l}(x)\vert^{2}&=
\sum^{l-1}_\zs{t=1}
\E_\zs{\p}\,\vert\tau_\zs{t,l}\vert^{2}
\le
\sum^{\q}_\zs{t=1}
\E_\zs{\p}\,\vert\tau_\zs{t,l}\vert^{2}
\\[2mm]
&=
\sigma^{4}_\zs{\p}
\sum_\zs{j,k}\,x_\zs{j}\,x_\zs{k}\,\Chi_\zs{\{\nu_\zs{j}=\nu_\zs{k}\}}
\bar{\psi}_\zs{k,l}\,
\psi_\zs{j,l}\,
\sum^{\q}_\zs{t=1}
\bar{\psi}_\zs{j,t}\,
\psi_\zs{k,t}\,
\,.
\end{align*}
Taking into account here the property \eqref{prd_scale}, we obtain that
$$
\E_\zs{\p}\vert\U_\zs{2,l}(x)\vert^{2}\le \q \sigma^{4}_\zs{\p} 
\sum_\zs{k}\,x^{2}_\zs{k}\vert\Upsilon_\zs{k}(\Delta_\zs{k})\vert^{4}
\le 2^{8}\x_\zs{*}^{2}\sigma^{4}_\zs{\p}\q^{-3}
\vert x\vert^{2}
$$
and, therefore,
$$
\E_\zs{\p}\,\vert \U_\zs{2}(x)\vert^{2}\le 2^{8}\x_\zs{*}^{2}\sigma^{4}_\zs{\p}\q^{-2}
\le
2^{8}\x_\zs{*}^{2}
\E_\zs{\p}\,\xi^{4}_\zs{1}
\q^{-2}
\,.
$$
From here it follows
\eqref{crl_12_0++}. Hence Proposition
\ref{Pr.sec:Prsm.3}.

\qed

\bigskip

\bigskip

\section{Proofs}\label{sec:Pr}

We will prove here most of the results of this paper.

\subsection{Proof of Theorem \ref{Th.sec:OI.1}} 

First, note that from 
\eqref{sec:Or.8} - \eqref{sec:Or.9}
it follows that
\begin{equation}\label{sec:Pr.1}
\Er_\zs{n}(\lambda) = J_n(\lambda) + 2 \sum_\zs{j\in\bbz^{d}}\, \lambda(j) \check{\theta}_\zs{j}+ \Vert S\Vert^2-\delta \wh{P}_n(\lambda)\,,
\end{equation}
where $\check{\theta}_\zs{j}=\wt{\theta}_\zs{j}-\Rep\,\bar{\theta}_\zs{j}\wh{\theta}_\zs{j}$. Using the definition of $\wt{\theta}_\zs{j}$ in \eqref{coeff_sigma}
we obtain that
\begin{align*}
\check{\theta}_\zs{j}&=
\Rep\,
\bar{\theta}_\zs{j}\, \zeta_\zs{j}
+
\vert \zeta_\zs{j}\vert^{2}
-\frac{\wh{\sigma}}{\q} \, \varpi_\zs{j}
=
\Rep\,
\bar{\theta}_\zs{j}\, \zeta_\zs{j}
+
\vert \b_\zs{j}\vert^{2}
\\[2mm]&
+\frac{2}{\sqrt{\q}}\,\Rep\,\bar{\b}_\zs{j} \eta_\zs{j}
+
\frac{1}{\q}\,\wt{\eta}_\zs{j}
+
\frac{\sigma_\zs{\p}-\wh{\sigma}}{\q} \, \varpi_\zs{j}
\,,
\end{align*}
where $\wt{\eta}_\zs{j}$ is defined in \eqref{func_U_++}.  Now we set 
\begin{align}\nonumber
\M(\lambda)& = \sum_\zs{j\in\bbz^{d}}\,\lambda(j)\Rep\,\bar{\theta}_\zs{j} \zeta_\zs{j}\,,
\qquad
\D_\zs{1}(\lambda)=\,
\frac{2}{\sqrt{\q}}\,
\sum_\zs{j\in\bbz^{d}}\,\lambda(j)
\,\Rep\,\bar{\b}_\zs{j} \eta_\zs{j}
\,,
\\[2mm]\label{sec:Pr.2}
\D_\zs{2}(\lambda)&=
\sum_\zs{j\in\bbz^{d}}\,\lambda(j)
\vert \b_\zs{j}\vert^{2}
\quad\mbox{and}\quad
\D(\lambda)=\D_\zs{1}(\lambda)+\D_\zs{2}(\lambda)
\,.
\end{align}
Using these functions,
we can rewrite \eqref{sec:Pr.1} as
\begin{align}\nonumber
\Er_\zs{n}(\lambda)  = &  J_n(\lambda) + 2 \frac{\sigma_\zs{\p}-  \wh{\sigma} }{\q}\,\check{L}(\lambda)+ 2 
\M(\lambda)+ 2 \D(\lambda)\\  \label{sec:Pr.3}
& +  2 \sqrt{P_\zs{n}(\lambda)} 
\frac{\check{\varkappa}(\lambda)\,\U(\e(\lambda))}{\sqrt{\sigma_\zs{\p} \q}}
 + \Vert S\Vert^2-\delta  P_n(\lambda),
\end{align}
where $\e(\lambda)=\lambda/\vert\lambda\vert$,  
$\check{\varkappa}(\lambda)=\vert\lambda\vert/\sqrt{\check{L}(\lambda^{2})}$ and the function $\check{L}(\cdot)$ is defined in \eqref{sec:Mo.2}. 
Let now $\lambda_\zs{0}= (\lambda_\zs{0}(j))_\zs{1\le j\le\,n}$ be a fixed sequence in $\Lambda$, $\wh{\lambda}$ be as in \eqref{sec:Mo.9}
and $\mu_\zs{0}=\wh{\lambda}-\lambda_\zs{0}$.
Substituting $\lambda_\zs{0}$ and $\wh{\lambda}$ in Equation \eqref{sec:Pr.3}, we obtain
\begin{align}\label{sec:Pr.4}
\Er_\zs{n}(\wh{\lambda})-\Er_\zs{n}(\lambda_\zs{0}) = & J(\wh{\lambda})-J(\lambda_\zs{0})+
2 \frac{\sigma_\zs{\p}-\wh{\sigma}}{\q}\,\check{L}(\mu_\zs{0})
+2 \M(\mu_\zs{0})
+2 \D(\mu_\zs{0})
\nonumber\\[2mm]
& + 2 \sqrt{P_\zs{n}(\wh{\lambda})} \frac{\check{\varkappa}(\wh{\lambda})\U(\wh e)}{\sqrt{\sigma_\zs{\p} \q}}-2 \sqrt{P_\zs{n}(\lambda_\zs{0})}
 \frac{\check{\varkappa}(\lambda_\zs{0})\U(e_0)}{\sqrt{\sigma_\zs{\p} \q}}\nonumber \\[2mm]
& -  \delta  \wh{P}_n(\wh{\lambda})+\delta \wh{P}_n(\lambda_\zs{0}),
\end{align}
where  $\wh{e} = e(\wh{\lambda})$ and $\e_\zs{0} = \e(\lambda_\zs{0})$. Note that, by \eqref{sec:Mo.2},
$$ 
|\check{L}(\mu_\zs{0})| \le\,\check{L}(\wh{\lambda}) + \check{L}(\lambda) \leq 2\vert\Lambda\vert_\zs{*}. 
$$
Using the inequality
\begin{equation}\label{sec:Pr.5}
2|ab| \leq \delta a^2 + \delta^{-1} b^2
\end{equation}
and taking into account that $P_\zs{n}(\lambda)> 0$ we obtain that
for any $\lambda\in\Lambda$ and any $0<\check{\delta}\le \delta$
$$
2 \sqrt{P_\zs{n}(\lambda)} \frac{\check{\varkappa}(\lambda)|\U(\e(\lambda))|}{\sqrt{\sigma_\zs{\p} \q}} \le\, \check{\delta} P_\zs{n}(\lambda) +
 \frac{\check{\varkappa}^{*}\U^{*}}{\check{\delta}\sigma_\zs{\p}\,\q}\,,
$$
where $\check{\varkappa}^{*}=\max_\zs{\lambda\in\Lambda}\check{\varkappa}^{2}(\lambda)$ and $\U^{*} = \sup_\zs{\lambda\in\Lambda} \U^{2}((\e(\lambda))$. 
Note here that for any $\lambda\in\Lambda$
$$
\vert
\wh{P}_\zs{n}(\lambda)
-
P_\zs{n}(\lambda)
\vert
\le
\frac{ |\wh{\sigma} -\sigma_\zs{\p}|}{\q}\,
\check{L}(\lambda^{2})
\le 
\frac{ |\wh{\sigma} -\sigma_\zs{\p}|}{\q}\,
\vert \Lambda\vert_\zs{*}
\,.
$$
So, taking into account  that $J(\wh{\lambda})\le J(\lambda_\zs{0})$, we get for any $0<\check{\delta}\le \delta<1$
that
\begin{align}\nonumber
\Er_\zs{n}(\hat \lambda)  \le &\,\Er_\zs{n}(\lambda_\zs{0})
+
6 \frac{\vert \wh{\sigma}-\sigma_\zs{\p}\vert}{\q}\,\vert \Lambda\vert_\zs{*} +2 \M(\mu_\zs{0})
+2 \D(\mu_\zs{0})
\\[2mm]\label{sec:Pr.6}
& 
+ \frac{2\check{\varkappa}^{*} \U^{*}}{\check{\delta}\sigma_\zs{\p}\q} 
+
(\check{\delta}-\delta)\,P_\zs{n}(\wh{\lambda})
+ 2  \delta P_\zs{n}(\lambda_\zs{0})\,.
\end{align}
To estimate the third term in the right side of this inequality we it represent for any 
$\mu\in\Lambda_\zs{1}=\Lambda-\lambda_\zs{0}=\left\{\lambda - \lambda_\zs{0}\,,\,\lambda\in \Lambda \right\}$ as
\begin{equation}
\label{sec:Mrepre__}
\M(\mu)=\M_\zs{1}(\mu)+\M_\zs{2}(\mu)
\,,
\end{equation}
where 
$$
\M_\zs{1}(\mu)=\frac{1}{\sqrt{\q}}\,\Rep\,
\sum_\zs{j\in\bbz^{d}}\,\mu(j)\,\bar{\theta}_\zs{j} \eta_\zs{j}
\quad\mbox{and}\quad
\M_\zs{2}(\mu)=
\sum_\zs{j\in\bbz^{d}}\,\mu(j)\Rep\,\bar{\theta}_\zs{j} \b_\zs{j}
$$
Moreover, for any family  $\upsilon=(\upsilon(j))_\zs{j\in\bbz^{d}}$ 
for which $\vert \upsilon\vert^{2}=\sum_\zs{j\in\bbz^{d}} \upsilon^{2}(j)
<\infty $
we set
\begin{equation}
\label{sec:S__projection}
S_\zs{\upsilon}(x) = \sum_\zs{j\in\bbz^{d}} \upsilon(j) \theta_\zs{j}  \Phi_\zs{j}(x)\,.
\end{equation}
Using 
Proposition \ref{Pr.sec:Prsm.2}
we obtain that
\begin{equation}\label{sec:Pr.8}
\E_\zs{\p} \vert\M_\zs{1}(\mu)\vert^{2} \le\,
 \frac{\check{\varpi}_\zs{\p}}{\q}\,\sum_\zs{j\in\bbz^{d}}\mu^{2}(j)\vert\theta_\zs{j}\vert^{2} 
 =\check{\varpi}_\zs{\p}\,
 \frac{\Vert\,S_\zs{\mu}\Vert^2}{\q}\,.
\end{equation}
To estimate this function for a random family $\mu=(\mu(j))_\zs{j\in\bbz^{d}}$ we set
\begin{equation}
\label{Z_setting_1}
Z^* = \sup_\zs{x \in \Lambda_1} 
\,\frac{\q \vert \M_\zs{1}(x)\vert^{2}}{\Vert\,S_\zs{x}\Vert^2}\,.
\end{equation}
So, through the inequality  \eqref{sec:Pr.5}, we get
\begin{equation*}\label{sec:Pr.10--}
2 |\M_\zs{1}(\mu)|\leq \delta \Vert\,S_\zs{\mu}\Vert^2 + \frac{Z^*}{\q\delta}.
\end{equation*}
It is clear that the last term  here can be estimated as
\begin{equation}\label{sec:Pr.9}
\E_\zs{\p} Z^* \leq \sum_\zs{x \in \Lambda_1} \frac{\q \E_\zs{\p} \vert\M_\zs{1}(x)\vert^{2}}{\Vert\,S_\zs{x}\Vert^2} 
\leq \sum_\zs{x \in \Lambda_1} \check{\varpi}_\zs{\p}= \check{\varpi}_\zs{\p}\check{\iota}\,,
\end{equation}
where $\check{\iota} = \mbox{card}(\Lambda)$. 
Using again the inequality \eqref{sec:Pr.5} we obtain that for any $x\in\Lambda_\zs{1}$
\begin{equation}\label{sec:Pr.13--1}
2 \vert\M_\zs{2}(x)\vert \leq \delta \Vert\,S_\zs{x}\Vert^2 + \frac{1}{\delta}
\sum_\zs{j\in\bbz^{d}}\,\vert x(j)\vert
\vert \b_\zs{j}\vert^{2}
\,
\le 
 \delta \Vert\,S_\zs{x}\Vert^2 + \frac{2 \vert \Lambda\vert_\zs{*}
\b_\zs{*}}{\delta}
\,,
\end{equation}
where
$
\b_\zs{*}=\sup_\zs{j\in\bbz^{d}}\vert\b_\zs{j}\vert^{2}$.
Thus,
\begin{equation}\label{sec:Pr.10}
2 |\M(\mu|\leq 2\delta \Vert\,S_\zs{\mu}\Vert^2 + \frac{Z^*}{\q\delta}
+\frac{2\vert \Lambda\vert_\zs{*}\b_\zs{*}}{\delta}
\,.
\end{equation}
Moreover, note that, for any $x\in\Lambda_1$,
\begin{align}\nonumber
\Vert\,S_\zs{x}\Vert^2-\Vert\wh{S}_x\Vert^2& = \sum_\zs{j\in\bbz^{d}}\, x^2(j) (\vert \theta_\zs{j}\vert^{2}-\vert\wh{\theta}_\zs{j}\vert^{2})\\[2mm] \label{sec:Pr.11}
& \le -2 \M_\zs{1}(x^{2})
-2 \M_\zs{2}(x^{2})\,,
\end{align}
where $x^{2} =(x^{2}(j))_\zs{j\in\bbz^{d}}$.
Taking into account that, for any $x \in \Lambda_1$ the components $|x(j)|\leq 1$,  we can estimate this term as
in \eqref{sec:Pr.8}, i.e.,
$$
\E_\zs{\p}\, \vert\M_\zs{1}(x^{2})\vert^{2} \leq \check{\varpi}_\zs{\p}\,
\frac{\Vert\,S_\zs{x}\Vert^2}{\q}\,.
$$
Similarly to the previous reasoning
we set
$$ 
Z^*_\zs{1} = \sup_\zs{x \varepsilon \Lambda_1} \frac{\q \vert\M_\zs{1}(x^{2})\vert^{2}}{\Vert\,S_\zs{x}\Vert^2}
$$
and we get
\begin{equation}\label{sec:Pr.12}
\E_\zs{\p}\, Z^*_1 \leq \check{\varpi}_\zs{\p}\,\check{\iota}\,.
\end{equation}
Using the same  arguments as in \eqref{sec:Pr.10}, we can derive
\begin{equation*}\label{sec:Pr.13}
2 \vert\M_\zs{1}(x^{2})\vert \leq \delta \Vert\,S_\zs{x}\Vert^2 + \frac{Z^*_1}{\q\delta}.
\end{equation*}
Similarly to \eqref{sec:Pr.13--1} we can estimate 
$$
2 \vert\M_\zs{2}(x^{2})\vert
\,\le 
 \delta \Vert\,S_\zs{x}\Vert^2 + \frac{2 \vert \Lambda\vert_\zs{*}
\b_\zs{*}}{\delta}
\,.
$$
From here and \eqref{sec:Pr.11} we get
\begin{equation}\label{sec:Pr.14}
\Vert\,S_\zs{x}\Vert^2 \leq \frac{\Vert\wh{S}_x\Vert^2}{1-2\delta} + \frac{Z^*_1}{\q\delta (1-2\delta)}
+ \frac{2 \vert \Lambda\vert_\zs{*}\b_\zs{*}}{\delta (1-2\delta)}
\end{equation}
for any $0<\delta<1$. Using this bound in \eqref{sec:Pr.10} yields
$$
2\vert \M(x)\vert \leq \frac{2\delta \Vert\wh{S}_x\Vert^2}{1-2\delta} + \frac{Z^*+Z^*_1}{\q \delta (1-2\delta)}
+
\frac{2 \vert \Lambda\vert_\zs{*}\b_\zs{*}}{1-2\delta}
\,.
$$
Taking into account that 
$$
\Vert\wh{S}_\zs{\mu_\zs{0}}\Vert^{2}= 
\Vert
\wh{S}_\zs{\wh{\lambda}} - 
\wh{S}_\zs{\lambda_\zs{0}}
\Vert^{2}
=
\Vert
(\wh{S}_\zs{\wh{\lambda}}-S) - 
(\wh{S}_\zs{\lambda_\zs{0}}-S)
\Vert^{2}
\le 2\,(\Er_\zs{n}(\wh{\lambda})+\Er_\zs{n}(\lambda_\zs{0}))
\,,
$$
we obtain
$$
2 \vert\M(\mu_\zs{0})\vert \leq \frac{2\delta(\Er_\zs{n}(\wh{\lambda})+\Er_\zs{n}(\lambda_\zs{0}))}{1-2\delta} + \frac{Z^*+Z^*_1}{\q \delta (1-2\delta)}
+
\frac{2 \vert \Lambda\vert_\zs{*}\b_\zs{*}}{1-2\delta}
\,.
$$
Let us estimate now the term $\D(\mu_\zs{0})$. Using the inequality \eqref{sec:Pr.5}
we obtain that for any $\lambda\in\Lambda$
 and $0<\check{\delta}<\delta<1$
\begin{align*}
\vert\D_\zs{1}(\lambda)\vert
&\le 
\frac{\check{\delta}}{\q}
\sum_\zs{j\in\bbz^{d}}\,\lambda^{2}(j)
\,\vert \eta_\zs{j}\vert^{2}
+
\frac{1}{\check{\delta}} 
\sum_\zs{j\in\bbz^{d}}\,\Chi_\zs{\{\lambda(j)\neq 0\}}
\,\vert\b_\zs{j} \vert^{2}\\[2mm]
&\le \check{\delta}P_\zs{n}(\lambda)
+
\frac{\check{\delta}}{\q}\U(\lambda^{2})
+
\frac{\vert\Lambda\vert_\zs{*}\b_\zs{*}}{\check{\delta}}
\,.
\end{align*}
Taking into account here that
\begin{align*}
\frac{\vert \U(\lambda^{2})\vert}{\q}
&\le  
\frac{\vert \lambda\vert\,\vert \U(\e(\lambda^{2}))\vert}{\q}
=\sqrt{P_\zs{n}(\lambda)}\,\check{\varkappa}(\lambda)\,
\frac{\vert\U(\e(\lambda^{2}))\vert}{\sqrt{\sigma_\zs{\p}\q}}\\[2mm]
&\le 
\check{\delta}\,P_\zs{n}(\lambda)
+
\,
\frac{\check{\varkappa}^{*}\U^{*}_\zs{1}}{\sigma_\zs{\p}\q\check{\delta}}
\,,
\end{align*}
where  $\U^{*}_\zs{1} = \sup_\zs{\lambda\in\Lambda} \U^{2}((\e(\lambda^{2}))$.
This implies that for any $\lambda\in\Lambda$
$$
\vert\D_\zs{1}(\lambda)\vert
\le 
2\check{\delta}\,P_\zs{n}(\lambda)
+
\,
\frac{\check{\varkappa}^{*}\U^{*}_\zs{1}}{\sigma_\zs{\p}\q}
+
\frac{\vert\Lambda\vert_\zs{*}\b_\zs{*}}{\check{\delta}}
$$
and, therefore,
$$
\vert\D_\zs{1}(\mu_\zs{0})\vert
\le 
2\check{\delta}\,P_\zs{n}(\wh{\lambda})
+
2\check{\delta}\,P_\zs{n}(\lambda_\zs{0})
\,+
\frac{2\check{\varkappa}^{*}\U^{*}_\zs{1}}{\sigma_\zs{\p}\q}
+
\frac{2\vert\Lambda\vert_\zs{*}\b_\zs{*}}{\check{\delta}}
\,.
$$
Moreover, similarly to the upper bound  \eqref{sec:Pr.13--1} we get
$$
\vert \D_\zs{2}(\mu_\zs{0})\vert \le \max_\zs{\lambda\in\Lambda}\D_\zs{2}(\lambda)+\D_\zs{2}(\lambda_\zs{0})
\le 
2 \vert \Lambda\vert_\zs{*}
\b_\zs{*}
\,.
$$
Finally, we obtain that
\begin{equation}
\label{Dmu0-Upp}
2\vert \D(\mu_\zs{0})\vert \le 
4\check{\delta}\,P_\zs{n}(\wh{\lambda})
+
4\check{\delta}\,P_\zs{n}(\lambda_\zs{0})
\,+
\frac{4\check{\varkappa}^{*}\U^{*}_\zs{1}}{\sigma_\zs{\p}\q}
+
\frac{8\vert\Lambda\vert_\zs{*}\b_\zs{*}}{\check{\delta}}
\,.
\end{equation}
So, using the upper bound \eqref{sec:Pr.6}, we obtain that
\begin{align*}\nonumber
\Er_\zs{n}(\hat \lambda)  \le &\,
\Er_\zs{n}(\lambda_\zs{0})
+
 \frac{2\delta(\Er_\zs{n}(\wh{\lambda})+\Er_\zs{n}(\lambda_\zs{0}))}{1-2\delta}
+
6 \frac{\vert \wh{\sigma}-\sigma_\zs{\p}\vert}{\q}\,\vert \Lambda\vert_\zs{*} \\[2mm]
& + \frac{Z^*+Z^*_1}{\q \delta (1-2\delta)}
+
\frac{8 \vert \Lambda\vert_\zs{*}\b_\zs{*}}{(1-2\delta)\check{\delta}}
\\[2mm]
& 
+ \frac{8\check{\varkappa}^{*} (\U^{*}+\U^{*}_\zs{1})}{\check{\delta}\sigma_\zs{\p}\q} 
+
(5\check{\delta}-\delta)\,P_\zs{n}(\wh{\lambda})
+ 6  \delta P_\zs{n}(\lambda_\zs{0})\,.
\end{align*}
So, choosing $\check{\delta}=\delta/5$, we obtain that

\begin{align*}\nonumber
\Er_\zs{n}(\hat \lambda)  \le &\,
 \frac{\Er_\zs{n}(\lambda_\zs{0})}{1-4\delta}
+
6 \frac{\vert \wh{\sigma}-\sigma_\zs{\p}\vert}{\q(1-4\delta)}\,\vert \Lambda\vert_\zs{*} 
 + \frac{Z^*+Z^*_1}{\q \delta (1-4\delta)}
+
\frac{40 \vert \Lambda\vert_\zs{*}\b_\zs{*}}{(1-4\delta) \delta}
\\[2mm]
& 
+ \frac{40\check{\varkappa}^{*} (\U^{*}+\U^{*}_\zs{1})}{\delta \sigma_\zs{\p}\q} 
+ \frac{6\delta}{1-4\delta}\, P_\zs{n}(\lambda_\zs{0})\,.
\end{align*}
In view of Proposition \ref{Pr.sec:Prsm.3} we estimate the expectation of the term $\U^{*}+\U^{*}_\zs{1}$  as
$$
\E_\zs{\p}\, (\U^{*}+\U^{*}_\zs{1}) \leq \sum_\zs{\lambda\in\Lambda}
\left(
\E_\zs{\p} \U^{2} (\e(\lambda))
+
\E_\zs{\p} \U^{2} (\e(\lambda^{2}))
\right) \leq 2\check{\iota} \c_\zs{*}\E_\zs{\p}\xi^{4}_\zs{1}\,.
$$
Taking into account that $0<\delta\le 1/8$, we get
\begin{align*}
 \cR(\wh{S}_*,S)  \le & \frac{ \cR(\wh{S}_\zs{\lambda_\zs{0}},S)}{1-4\delta} 
+  \frac{12\vert\Lambda\vert_\zs{*}\,\E_\zs{\p}\,|\wh{\sigma} -\sigma_\zs{\p}|}{\q}
+ \frac{4\check{\varpi}_\zs{\p} \check{\iota}}{\q \delta}+
\frac{80 \vert \Lambda\vert_\zs{*}\b_\zs{*}}{\delta}
\\[4mm]
& 
+ \frac{80\check{\varkappa}^{*} \check{\iota} \c_\zs{*}\E_\zs{\p}\xi^{4}_\zs{1}}{\delta \sigma_\zs{\p}\q} 
  + \frac{2\delta}{(1-4\delta)} P_\zs{n}(\lambda_\zs{0}).
\end{align*}
Using the upper bound for $ P_n(\lambda_\zs{0})$ in Lemma~\ref{Le.sec:A.1-06-11-01}, one obtains 
\eqref{Th.sec:OI.1}, that finishes the proof. \fdem

\subsection{Proof of Proposition \ref{Pr.sec:Si.1}}

 We use here the same method as in \cite{KonevPergamenshchikov2009a}.
First of all note that the definition
\eqref{equa2.8} implies that
\begin{equation}\label{sec:Mo.1-0-1-04}
\wh{\sigma}=
\frac{1}{\wt{\q}}
\,
\sum_\zs{j\in\cT_\zs{n}}\,\vert \theta_\zs{j}\vert^{2}
+
\frac{2}{\wt{\q}}
\,
\M_\zs{n}
+
\frac{1}{\wt{\q}}
\,
\sum_\zs{j\in\cT_\zs{n}}\,
\vert\zeta_\zs{j}\vert^{2}
\,,
\end{equation}
where $\M_\zs{n}=
\,
\Rep
\sum_\zs{j\in\cT_\zs{n}}\,
\bar{\theta}_\zs{j}\,
\zeta_\zs{j}$.
In Lemma \ref{Le.sec:A.FC++} we show that
\begin{equation}\label{sec:App.2-1-04} 
\sum_\zs{j\in \cQ_\zs{n}}\,\vert \theta_\zs{j}\vert^{2}
\le 
\frac{\x^{2d}_\zs{*}}{\pi^{2d}\,\m_\zs{n}^{d/2}}
\,
\check{\tau}_\zs{d}(S)
\,,
\end{equation}
where $\cQ_\zs{n}=\{[\sqrt{\m_\zs{n}}]+1,\ldots,\}^{d}$.
To estimate the second term in \eqref{sec:Mo.1-0-1-04}
we represent it as
$$
\M_\zs{n}=\M_\zs{1,n}+\M_\zs{2,n}\,,
$$
where
$$
\M_\zs{1,n}=\Rep\sum_\zs{j\in\cT_\zs{n}}\,
\bar{\theta}_\zs{j}\,
\b_\zs{j}
\quad\mbox{and}\quad
\M_\zs{2,n}=\frac{1}{\sqrt{\q}}
\Rep
\sum_\zs{j\in\cT_\zs{n}}\,
\bar{\theta}_\zs{j}\,
\eta_\zs{j}
\,.
$$
To estimate $\M_\zs{1,n}$
note that
$$
2\vert \M_\zs{1,n}\vert
\le
\,
\sum_\zs{j\in\cT_\zs{n}}\,
\vert\theta_\zs{j}\vert^{2}
+
\,
\sum_\zs{j\in\cT_\zs{n}}\,
\vert\b_\zs{j}\vert^{2}\,
\le\,
\frac{\x^{2d}_\zs{*}}{\pi^{2d}\,\m_\zs{n}^{d/2}}
\,
\check{\tau}_\zs{d}(S)
+
\sum_\zs{j\in\cT_\zs{n}}\,
\vert\b_\zs{j}\vert^{2}\,
\,.
$$
To estimate the last term in this inequality 
we use the coefficient
 $\wt{\b}$ defined in
\eqref{psi-qq-p} and the fact that $\q\ge \m^{d}_\zs{n}$, i.e.
\begin{equation}
\label{bound-b-n}
\sum_\zs{j\in\cT_\zs{n}}\,
\vert\b_\zs{j}\vert^{2}
\le \wt{\b} \frac{\m^{d}_\zs{n}}{\q^{2}}
\le  \frac{\wt{\b}}{\q}
\,.
\end{equation}
Therefore,
$$
2\vert \M_\zs{1,n}\vert
\le
\frac{\x^{2d}_\zs{*}}{\pi^{2d}\,\m_\zs{n}^{d/2}}
\,
\check{\tau}_\zs{d}(S)
+
\frac{\wt{\b}}{\q}
\,.
$$
Moreover,
the term $\M_\zs{2,n}$
can be estimated through
Proposition \ref{Pr.sec:Prsm.2}
  as
$$
\E_\zs{\p}\,\vert\M_\zs{2,n}\vert^{2}\le \frac{\check{\varpi}_\zs{\p}}{\q}\,
\sum_\zs{j\in\cT_\zs{n}}\,\vert\theta\vert^{2}_\zs{j}
\le 
\frac{\check{\varpi}_\zs{\p}\check{\tau}_\zs{d}(S)\x^{2d}_\zs{*}}{\pi^{2d}\,\m_\zs{n}^{d/2}\q}
\,,
$$
while the absolute value of this term  can be estimated  as
$$
2
\E_\zs{\p}\,
\vert \M_\zs{2,n}\vert
\le 2
\frac{\sqrt{\check{\varpi}_\zs{\p}\check{\tau}_\zs{d}(S)}\x^{d}_\zs{*}}{\pi^{d}\,\m_\zs{n}^{d/4}\sqrt{\q}}
\,
\le 
\,
\frac{\check{\varpi}_\zs{\p}\check{\tau}_\zs{d}(S)\x^{2d}_\zs{*}}{\pi^{2d}\,\m_\zs{n}^{d/2}}
+
\frac{1}{\q}
\,.
$$
Therefore,
$$
2
\frac{\E_\zs{\p}\,
\vert \M_\zs{n}\vert}{\wt{\q}}
\le 
\frac{(1+\check{\varpi}_\zs{\p})\check{\tau}_\zs{d}(S)\x^{2d}_\zs{*}}{\pi^{2d}\,\m_\zs{n}^{d/2}\wt{\q}}
+
\frac{1+\wt{\b}}{\q\wt{\q}}
\,.
$$
We can represent the last term in \eqref{sec:Mo.1-0-1-04} as
\begin{equation}
\label{zeta++2}
\frac{1}{\wt{\q}}
\,
\sum_\zs{j\in\cT_\zs{n}}\,
\vert\zeta_\zs{j}\vert^{2}
=
\frac{1}{\wt{\q}\q}
\sum_\zs{j\in\cT_\zs{n}}\,
\vert\eta_\zs{j}\vert^{2}
+
\frac{1}{\wt{\q}}
\,
\sum_\zs{j\in\cT_\zs{n}}\,
\vert\b_\zs{j}\vert^{2}
+\frac{2}{\wt{\q}\sqrt{\q}}
\Rep 
\sum_\zs{j\in\cT_\zs{n}}\,
\eta_\zs{j}\bar{\b}_\zs{j}
\,.
\end{equation}
Moreover,  using the definition \eqref{func_U_++} we obtain
$$
\frac{1}{\q}
\sum_\zs{j\in\cT_\zs{n}}\,\vert\eta_\zs{j}\vert^2
=
\frac{\sigma_\zs{\p}\sum_\zs{j\in \cT_\zs{n}}\varpi_\zs{j}}{\q}
+\frac{\m_\zs{n}^{d/2}\U(\check{\x})}{\q}
=
\,
\sigma_\zs{\p}\wt{\q}
+\frac{\m_\zs{n}^{d/2}\U(\check{\x})}{\q}
$$
with
$\check{\x}_\zs{j}=\Chi_\zs{\{j\in \cT_\zs{n}\}}/\m_\zs{n}^{d/2}$.
Therefore, from Proposition  \ref{Pr.sec:Prsm.3} it follows that
$$
\E_\zs{\p}
\left\vert
\frac{1}{\wt{\q}\q}
\sum_\zs{j=\cT_\zs{n}}\,\vert\eta_\zs{j}\vert^2
-\sigma_\zs{\p}
\right\vert
\le 
\,
\frac{\m_\zs{n}^{d/2}\,\sqrt{\c_\zs{*}\,\E_\zs{\p}\xi^{4}_\zs{1}}}{\q\wt{\q}}
\,.
$$
Moreover, using here
\eqref{infvarpi--}, we obtain that for $\m_\zs{n}\ge 4$
$$
\E_\zs{\p}
\left\vert
\frac{1}{\wt{\q}\q}
\sum_\zs{j=\cT_\zs{n}}\,\vert\eta_\zs{j}\vert^2
-\sigma_\zs{\p}
\right\vert
\le 
\,
\frac{\m_\zs{n}^{d/2}\,\sqrt{\c_\zs{*}\,\E_\zs{\p}\xi^{4}_\zs{1}}}{\sqrt{\varpi_\zs{*}}(\m_\zs{n}-\sqrt{\m_\zs{n}})^{d}}
\le 
\frac{2^{d}\,\sqrt{\c_\zs{*}\,\E_\zs{\p}\xi^{4}_\zs{1}}}{\sqrt{\varpi_\zs{*}}\,\m^{d/2}_\zs{n}}
\,.
$$
To estimate the last term in \eqref{zeta++2} we use this bound \eqref{bound-b-n} and again Prposition \ref{Pr.sec:Prsm.2}, i.e.
$$
\E\,\left\vert
\sum_\zs{j\in\cT_\zs{n}}\,
\eta_\zs{j}\bar{\b}_\zs{j} \right\vert^{2}
\le 
\,
\check{\varpi}_\zs{\p}
\,
\sum_\zs{j\in\cT_\zs{n}}\,
\vert\b_\zs{j}\vert^{2}
\le
\,
 \frac{\check{\varpi}_\zs{\p}\wt{\b}}{\q}
\,.
$$
Thus,
\begin{align*}
\E_\zs{\p}
\left\vert
\frac{1}{\wt{\q}}
\,
\sum_\zs{j\in\cT_\zs{n}}\,
\vert\zeta_\zs{j}\vert^{2}
-
\sigma_\zs{\p}
\right\vert
&\le
\frac{2^{d}\,\sqrt{\c_\zs{*}\,\E_\zs{\p}\xi^{4}_\zs{1}}}{\sqrt{\varpi_\zs{*}}\,\m^{d/2}_\zs{n}}
+
 \frac{\wt{\b}}{\wt{\q}\q}
 +
 \frac{2\sqrt{\check{\varpi}_\zs{\p}\wt{\b}}}{\wt{\q}\,\q}
\\[3mm]
&\le 
\frac{2^{d}\,\sqrt{\c_\zs{*}\,\E_\zs{\p}\xi^{4}_\zs{1}}}{\sqrt{\varpi_\zs{*}}\,\m^{d/2}_\zs{n}}
+
 \frac{2\wt{\b}+\check{\varpi}_\zs{\p}}{\wt{\q}\q}
 \,.
\end{align*}
It should be noted also that  for $\q\le \q_\zs{*}$
we can estimate $\wt{\q}$ from below as
$$
\wt{\q}=
\frac{\sum_\zs{j\in \cT_\zs{n}}\varpi_\zs{j}}{\q}
\ge \frac{\sqrt{\varpi_\zs{*}}}{2^{d}}\frac{\m^{d}_\zs{n}}{\q}
\ge 
\frac{\sqrt{\varpi_\zs{*}}}{2^{d}}\frac{1}{\q_\zs{*}}
>0\,.
$$
Moreover, using Lemma \ref{Le.sec.RdnTrs_1}
we can estimate directly $\b_\zs{j}$ as
$$
\sup_\zs{j\in\bbz^{d}}
\vert\b_\zs{j}\vert
\le \frac{16\check{\r}}{(2\x_\zs{*})^{d/2-2}}\,\frac{1}{\q}
\,.
$$
From here
we obtain the bound \eqref{sec:Si.3} and hence the desired result.
\qed

\bigskip

{\bf Acknowledgments.}  
The last author is partially supported 
 by RFBR Grant 16-01-00121, 
 by the  Ministry of Education and Science of the Russian Federation  in the framework of the research project 
 No 2.3208.2017/4.6,
  by the Russian Federal Professor program  
 (Project No 1.472.2016/1.4, Ministry of Education and Science of the Russian Federation)
and 
"The Tomsk State University competitiveness improvement programme"\ Grant 8.1.18.2018.

\bigskip

\renewcommand{\theequation}{A.\arabic{equation}}
\renewcommand{\thetheorem}{A.\arabic{theorem}}
\renewcommand{\thesubsection}{A.\arabic{subsection}}
\section{Appendix}\label{sec:A}
\setcounter{equation}{0}
\setcounter{theorem}{0}

\subsection{Properties of the Radon transformation} \label{Prpes_1}

First, we recall some basic definitions valid for $\bbr^{d}\to \bbr$ functions 
$f$ belonging to the Schwartz space 
 and which can be found, for 
instance, in 
\cite{Natterer2005}. The Fourier transform of $f$ is given by
\begin{equation}\label{equa1.1}
 T(f)(z) = 
 \frac{1}{(2 \, \pi)^{d/2}} \int_{\R^d} f(x) \, e^{-i \, x \cdot z} \, dx
\end{equation}
and its inverse Fourier transform by 
\begin{equation}\label{equa1.2}
 T^{-1}(f)(\eta) = 
 \frac{1}{(2 \, \pi)^{d/2}} \int_{\R^d} f(y) \, e^{i \, y \cdot \eta} \, dy
\end{equation}
where $\cdot$ denotes the inner product in $\R^d$. 

For any $\nu$ in the unit sphere $S$ and any $\varsigma \in \R$, the Radon 
transform of $f$ is defined by 
\begin{equation}\label{equa1.3}
 R(f)(\nu,\varsigma) = \int_{\nu^\perp} f(\varsigma \, \nu + y) \, dy \, , 
\end{equation}
where $\nu^\perp$ is the subspace orthogonal to $\nu$ in $\R^d$. 
Setting 
\begin{equation}\label{equa1.4}
 A_\nu(f)(\varsigma) = R(f)(\nu, \varsigma) \, , 
\end{equation}
it is easily seen that 
\begin{equation}\label{equa1.5}
 T(f)(\varsigma \, \nu) = T \circ A_\nu(f)(\varsigma) \, .  
\end{equation}


\noindent Taking into account that 
for $\vert\varsigma\vert\ge x_\zs{*}$ and for $y\in \nu^\perp$ the norm
$\|\varsigma \, \nu + y\|^2 = \|\varsigma\|^2 + \|y\|^2 \geq N^2$, we obtain that
\begin{equation}\label{FinRadonTr_1}
R(f)(\nu,\varsigma) = 0
\quad\mbox{for any}\quad 
\vert\varsigma\vert \geq N
\,. 
\end{equation}

The following lemma gives a Lipschitzian property of the Radon transform. 

\begin{lemma}\label{Le.sec.RdnTrs_1}
 Let $f$ be a Lipschitzian function from $\R^d$ into $\R$ with Lipschitz constant $\check{\r}$ and with 
 compact support a centered ball of radius $\x_\zs{*}$. 
 
 Then its Radon transform is Lipschitzian with Lipschitz constant $\check{\r} \, \x_\zs{*}^{d-1}$; more precisely, 
 for any $\nu \in S$ and any $(s_1,s_2) \in \R^2$, 
 \begin{equation}\label{equa3.11}
  |R(f)(\nu,s_1) - R(f)(\nu,s_2)| \leq \check{\r} \, \x_\zs{*}^{d-1} \, |s_1 - s_2| \, . 
 \end{equation}
\end{lemma}

\noindent {\sc Proof}
Let $\nu \in S$ and $(s_1,s_2) \in \R^2$. By definition of the Radon transform, we have 
\begin{eqnarray*}
 |R(f)(\nu,s_1) - R(f)(\nu,s_2)| 
 &=& 
 \left| \int_{\R^{d-1}} (f(\nu \, s_1 + y) - f(\nu \, s_2 + y)) \, dy\right| \\
 &=& 
 \left| \int_{|y|\leq N} (f(\nu \, s_1 + y) - f(\nu \, s_2 + y)) \, dy\right| \\
 &\leq&
 \int_{|y|\leq N} \left| f(\nu \, s_1 + y) - f(\nu \, s_2 + y)\right| dy \\
 &\leq&
 \check{\r} \, \x_\zs{*}^{d-1} \, |s_1 - s_2| \, . 
\end{eqnarray*}
Hence Lemma \ref{Le.sec.RdnTrs_1}. \qed

\subsection{Property of the Fourier coefficients}

\begin{lemma}\label{Le.sec:A.FC++}
Assume that the partial derivative $\partial^{d}/\partial x_\zs{1}\dots\partial x_\zs{d}$ of $S$
is continuos. Then the equality \eqref{sec:App.2-1-04} holds true.
\end{lemma}
\proof Integrating by parts we  obtain that
$$
 \theta_\zs{j} = 
 \frac{\x^{d}_\zs{*}}{i^{d}\pi^{d}\prod^{d}_\zs{l=1}j_\zs{l}}
 \int_{[-\x_\zs{*},\x_\zs{*}]^d} \frac{\partial^{d} }{\partial x_\zs{1}\dots\partial x_\zs{d}}S(z) \, \Phi_\zs{j}(z) \, \d z
 \,.
$$
So, applying here the  Bunyakovsky - Cauchy - Swarths inequality we obtain the upper bound \eqref{sec:App.2-1-04}.
Hence lemma \ref{Le.sec:A.FC++}.
\fdem

\subsection{Property of the penalty term}

\begin{lemma}\label{Le.sec:A.1-06-11-01}
For any $\lambda \in \Lambda$,
$$ 
P_\zs{n}(\lambda) \leq \E_\zs{\p} \Er_\zs{n}(\lambda)
$$
where the coefficient $P_\zs{n}(\lambda)$ was defined in \eqref{sec:Pr.2}.
\end{lemma}
\proof
 By the definition of $\Er_\zs{n}(\lambda)$
  one has
\begin{align*}
\E_\zs{\p}\Er_\zs{n}(\lambda)
&= \sum_\zs{j\in\bbz^{d}}\, \E_\zs{\p}\left\vert(\lambda(j)-1) \theta_\zs{j}+ \lambda(j)\zeta_\zs{j} \right\vert^2\\[2mm]
&=\sum_\zs{j\in\bbz^{d}}\, \E_\zs{\p}\left\vert(\lambda(j)-1) \theta_\zs{j}+ \lambda(j)\b_\zs{j}+\lambda(j)\frac{1}{\sqrt{\q}}\eta_\zs{j} \right\vert^2
\\[2mm]
&\ge 
 \frac{1}{\q}
\sum_\zs{j\in\bbz^{d}} \lambda^{2}(j)\E_\zs{\p}\vert \eta_\zs{j}\vert^{2}
=P_\zs{n}(\lambda)\,. 
\end{align*}                                                                              
Hence lemma \ref{Le.sec:A.1-06-11-01}.
\fdem

\bigskip


\end{document}